\documentclass[11pt]{article}
 \usepackage[top=1in,bottom=1in,left=1.5in,right=1.5in]{geometry}
  \usepackage{amsmath,amssymb}
  \usepackage{latexsym}
  \usepackage{graphicx,color}

  \newcommand{\C}{\mathbb{C}}

  \newcommand{\R}{\mathbb{R}}

  \newcommand{\U}{\mathbf{U}}
  \newcommand{\uu}{\mathbf{u}}

  \newcommand{\x}{\mathbf{x}}

  \newcommand{\lam}{\mbox{\boldmath{$\lambda$}}}

  \newcommand{\cF}{\mathcal{F}}

  \newcommand{\cP}{\mathcal{P}}

  \def\diag{\mathop{{\rm diag}}\nolimits}
  \newcommand{\hs}{\hspace*{\parindent}}
  \newcommand{\proof}{\hs \textbf{Proof.\ }}
  
  \newcommand{\tr}{\mathop{\mathrm{tr}}\nolimits}

  \newcommand{\trans}{^\top}
  \newcommand{\qed}{\hspace*{\fill} $\Box$\\}

  \newcommand{\rH}{\mathrm{H}}

  \newtheorem{theo}{\bfseries \hs Theorem}[section]
  
  \newtheorem{prop}[theo]{\bfseries \hs Proposition}
  \newtheorem{lemma}[theo]{\bfseries \hs Lemma}

  \numberwithin{equation}{section} 

 \renewcommand{\span}{\mathrm{span}}

\begin{document}

 \title{Equality in Wielandt's eigenvalue inequality}

 \author
 {Shmuel Friedland\footnote{Department of Mathematics, Statistics and Computer Science,
 University of Illinois at Chicago,
 Chicago, Illinois 60607-7045, USA,
 \emph{email}: friedlan@uic.edu. This work was supported by NSF grant DMS-1216393}}

 \date{March 22, 2015}
 \maketitle

 \begin{abstract}
 In  this paper we  give necessary and sufficient conditions for the equality case in Wielandt's eigenvalue inequality.
 \end{abstract}

 \noindent \emph{Keywords}: Lidski's theorem, Wielandt's eigenvalue inequality.

 \noindent {\bf 2010 Mathematics Subject Classification}: 15A18, 15A22, 15A42, 15B57.
 \section{Introduction}\label{sec:intro}
 For a positive integer $n$ let $[n]:=\{1,2,\ldots,n\}$. 
 Denote by $\rH_n$ the real space on $n\times n$ hermitian matrices.  For $A\in\rH_n$ let $\lambda_1(A)\ge \cdots\ge \lambda_n(A)$ be the $n$ eigenvalues
 of $A$, counted with their multiplicities.  Let $\lam(A)=(\lambda_1(A),\ldots,\lambda_n(A))\trans$ and let $\tr A=\sum_{i=1}^n \lambda_i(A)$ be the trace of $A$.  
 Denote by $\cP_n\subset \R^{n\times n}$ the group of permutation
 matrices.  In a short note \cite{Lid50} Lidskii announced the following result:  Let $A,B\in\rH_n$.  Then $\lam(A+B)-\lam(A)$ is in the convex hull spanned by $P\lam(B)$,
 where $P\in\cP_n$.  This fact is equivalent to the result that $\lam(A+B)-\lam(A)=O\lam( B)$ for some doubly stochastic matrix $O$.  That is, $\lam(A+B)-\lam(A)$ is majorized 
 by $\lam(B)$ \cite{HLP52}.  Since $\tr(A+B)=\tr A+\tr B$ the result of Lidskii is equivalent to the inequalities
 \begin{equation}\label{wielandt1}
 \sum_{j=1}^k \lambda_{i_j}(A+B)\le\sum_{j=1}^k \lambda_{i_j}(A) + \sum_{j=1}^k \lambda_j(B),
 \end{equation}
 for each $k\in [n-1]$ and distinct integers $i_1,\ldots,i_k$ in $[n]$.  This inequality was proved by Wielandt \cite{Wie56} by using max-min characterization of 
 $\sum_{j=1}^k \lambda_{i_j}(A)$.  

 The aim of this note is to give necessary and sufficient conditions for the equality 
 \begin{equation}\label{wieleq}
 \sum_{j=1}^k \lambda_{i_j}(A+B)=\sum_{j=1}^k \lambda_{i_j}(A) + \sum_{j=1}^k \lambda_j(B),
 \end{equation}
 for given integers $1\le i_1<\cdots <i_k \le n$.
 We also give a simple proof of the inequality \eqref{wielandt1} using a variation formula for the eigenvalues of the pencil $A(t)=A+tB$ for $t\in[0,1]$.

 We now summarize briefly the contents of the paper.  In \S2 we state preliminary results on hermitian matrices and pencils.
 In \S3 we state and prove the main result of this paper: Theorem \ref{Wielineqeq}, which gives necessary and sufficient conditions for \eqref{wieleq}.
 In the last section we comment on the main result of the paper. 

 \section{Preliminary results on hermitian pencils}\label{sec:wielineq}
 Fix $A,B\in\rH_n$.  Then $A(z):=A+z B, z\in\C$ is called a \emph{hermitian pencil}.  The $n$ eigenvalues of $A(z)$ are algebraic functions satisfying
 the equation $\det (\alpha I_n -A(z))=0$.  These eigenvalues are multivalued functions on $\C$, each one with $n$ branches at most,
 which are locally analytic except at a finite number of points $Z\subset \C$.  Furthermore, at each $z\in \C\setminus Z$ $A(z)$ has exactly $K=K(z)$ distinct eigenvalues
 $\gamma_1(z),\ldots,\gamma_K(z)$,
 and each eigenvalue $\gamma_i(z)$ has a fixed mulitplicity $M_i=M_i(z)$ for $i\in K$.  Moreover,  $|Z|\le n(n-1)$.  See \cite{MoF80, Fri81}.

 For $t\in\R$ the matrix $A(t)$ is hermitian.  We arrange  its eigenvalues in a decreasing order
 \begin{equation}\label{deflmabt}
 \lambda_1(t)\ge \cdots\ge\lambda_n(t), \quad t\in\R.
 \end{equation}
 Hence $\lambda_1(t),\ldots,\lambda_n(t)$ satisfy the equation $\det (\lambda I_n - A(t))=0$, and they are analytic on $\R\setminus Z$.
 Furthermore
 \begin{eqnarray}\label{lambdtmultip}
 &&\lambda_1(t)=\cdots=\lambda_{M_1}(t)>\lambda_{M_1+1}(t)=\cdots=\lambda_{M_1+M_2}(t)>\cdots>\\
 &&\lambda_{M_1+\ldots+M_{K-1}+1}(t)=\cdots=\lambda_{M_1+\ldots+M_K}(t),\; n=M_1+\cdots+M_K, \;t\in\R\setminus Z.\notag
 \end{eqnarray}
 
 Note 
 \begin{eqnarray}\label{defintIj}
 &&\R\setminus Z=\cup_{i=1}^N I_j, \;I_j=(a_{j-1}, a_j), j\in[N], \\
 &&-\infty=a_0<a_1<\cdots<a_{N-1}<a_N=\infty.\notag
 \end{eqnarray}
 
 We now recall a well known perturbation formula for eigenvalues of $A(z)$ at $z=0$:
 \begin{lemma}\label{pertform}  Assume that $A\in\rH_n$.  Suppose furthermore that $A$ has exactly $l\in[n]$ distinct eigenvalues of multiplicities  $n_1,\ldots,n_l\in[n]$:
 \begin{eqnarray}\label{mulipeigA}
 &&\lambda_1(A)=\cdots=\lambda_{m_1}(A)>\lambda_{m_1+1}(A) =\cdots=\lambda_{m_2}(A)>\\
 &&\cdots >\lambda_{m_{l-1}+1}(A)=\cdots =\lambda_{m_l}(A), \notag\\
 && m_0=0, \quad m_j=n_1+\cdots +n_j \textrm{ for } j\in [l].\label{defmj}
 \end{eqnarray}
 Assume that $B\in\rH_n$. Then it is possible to arrange the eigenvalues of the pencil $A(z)=A+zB$ for $|z|<r$, where $r$ is small, as $\alpha_1(z),\ldots,\alpha_n(z)$ such that 
 \begin{equation}\label{pertformalal}
 \alpha_j(z)=\lambda_j(A)+ z(\nu_j(A,B) +o(|z|)), \quad j\in[n].
 \end{equation}
 Assume that 
 \begin{equation}\label{eigvectA}
 A\uu_i=\lambda_i(A)\uu_i, \;\uu_i\in\C^n, \quad \uu_i^*\uu_j=\delta_{ij}, \textrm{ for } i,j\in [n].
 \end{equation}
 Then $\nu_{m_{i-1}+j}(A,B), j\in[n_i]$ are the eigenvalues of the hermitian matrix $[\uu_j^* B\uu_k]_{j=k=1}^{n_i}$ arranged in the decreasing order for $i\in[l]$.
 In particular, it is possible to choose an orthonormal system of eigenvectors of $A$ satisfying \eqref{eigvectA} such that $\nu_j(A,B)=\uu_j^*B\uu_j$
 for $j\in [n]$. 
 \end{lemma}
 The proof of this lemma follows from a well known perturbation formula for an eigenvalue $\alpha(x)$ of $F+xG$, where $F,G\in\C^{n\times n}$ and $\alpha(0)$ 
 is a geometrically simple eigenvalue of $F$ \cite{Kat82}.  This perturbation formula is elementary \cite{Fri78}.  
 This lemma  is also a simple consequence of Rellich's theorem \cite{Rel69}.   
 We now summarize the above result in the following known theorem:
 \begin{theo}\label{difforeigpen}  Let $A,B\in\rH_n$ and let $A(z)=A+zB$ for $z\in\C$.  Then there exists a finite set $Z\subset \C$, possibly empty, of cardinality at most 
 $n(n-1)$, such that the eigenvalues of $A(z)$ are multivalued analytic functions on $\C\setminus Z$.  The number of distinct eigenvalues of $A(z)$ is $K$, and
 the eigenvalue $\gamma_i(z)$ is of multiplicity $M_i$ for $z\in \C\setminus Z$ for $i\in [K]$.

 For $t\in\R$ arrange the eigenvalues of $A(t)$ as in \eqref{deflmabt}.   Assume that 
 intervals $I_j=(a_{j-1},a_{j})$ for $j\in[N]$ are given by \eqref{defintIj}.  Then \eqref{lambdtmultip} holds.
 For each $t\in\R$ there exists a choice of orthonormal eigenvectors of $A(t)$
 \begin{equation}\label{eigenvecAt}
 A(t)\uu_i(t)=\lambda_i(t)\uu_i(t), \quad \uu_i(t)^*\uu_k(t)=\delta_{ik} \textrm{ for } i,k\in [n],
 \end{equation}
 such that the following conditions hold:
 \begin{eqnarray}\label{difforeiglami}
 &&\lambda_i'(t)=\uu_i(t)^*B\uu_i(t), \quad i\in [n], \textrm{ for } t\in I_j,\\
 &&\lambda_i'(a_j^+)=\uu_i(a_j)^* B\uu_j(a_j), \quad i\in [n], \textrm { for } j\in [N-1]. \label{difforeiglami+}
 \end{eqnarray}
 Furthermore $\lambda_i'(t)$ is continuous from the right and from the left at $t=a_j$ for each $i\in[n]$ and $j\in [N-1]$. 
 In particular
 \begin{equation}\label{derlambi0+}
 \lambda_i'(0^+)=\nu_i(A,B), \quad i\in[n],
 \end{equation}
 where $\nu_1(A,B),\ldots,\nu_n(A,B)$ are defined in Lemma \ref{pertform}.

 \end{theo}

 Rellich's theorem \cite{Rel69} states that there exists a connected open set $\Omega\subset \C$ containing $\R$, such that $A(z)$
 has $n$ analytic eigenvalues $\alpha_1(z),\ldots,\alpha_n(z)$ and the corresponding analytic eigenvectors $\uu_1(z),\ldots,\uu_n(z)$ in $\Omega$.
 Furthermore, $\alpha_1(t),\ldots,\alpha_n(t)$ are real and $\uu_1(t),\ldots,\uu_n(t)$ are orthonormal for $t\in\R$. 
 
 Recall the Ky Fan charaterization of the sum of the first $k$-eigenvalues of $A\in\rH_n$ \cite{Fan49}.  Let $\cF_{k,n}$ be the set of all $k$ orthonormal
 vectors $\{\x_1,\ldots,\x_k\}$ in $\C^n$.  That is, $\x_i^*\x_j=\delta_{ij}$ for $i,j\in[k]$.  Then
 \begin{equation}\label{kyfanchar}
 \sum_{i=1}^k \lambda_i(A)=\max_{\{\x_1,\ldots,\x_k\}\in \cF_{k,n}}\sum_{i=1}^k \x_i^* A\x_i, \quad A\in\rH_n.
 \end{equation}
 Note that for $k=n$ we have the equality $\tr A=\sum_{i=1}^n \x_i^* A\x_i$ for each $\{\x_1,\ldots,\x_n\}\in\cF_{n,n}$.
 Equality in \eqref{kyfanchar} holds if and only if $\span(\x_1,\ldots,\x_k)$ is an invariant subspace of $A$ corresponding to the first $k$ eigenvalues of $A$.

 \section{A characterization of the equality case}\label{sec:mainthe}
 \begin{theo}\label{Wielineqeq}  Let $A,B\in\rH_n$ and $k\in [n-1]$.   Assume that $1\le i_1<\cdots<i_k\le n$.  Then \eqref{wielandt1} holds.
 Equality \eqref{wieleq} holds if and only if the following conditions are satisfied:
 There exist $r$ invariant subspaces $\U_1,\ldots,\U_r\subset\C^n$ of $A$ and $B$ such that
 each $\U_l$ is spanned by $k$-orthonormal vectors of $B$ corresponding to $\lambda_1(B),\ldots,\lambda_k(B)$.
 Let  $\mu_{1,l}(t)\ge \cdots\ge\mu_{k,l}(t)$ be the eigenvalues of the restriction of $A(t)$ to $\U_l$ for $l=1,\ldots,r$. 
 Then there exist $b_0=0<b_1<\cdots<b_{r-1}<b_r=1$ with the following properties: For each $l\in [r]$ and $t\in [b_{l-1},b_l]$
 $\mu_{j,l}(t)=\lambda_{i_j}(t)$ for $j=1,\ldots,k$. 
 \end{theo}
 \proof Let $I_l$ be an interval as in Theorem \ref{difforeigpen}.
 Let $\phi(t)=\sum_{j=1}^k \lambda_{i_j}(t)$ for $t\in\R$.  Then for $t\in I_l$ one has:
 \begin{equation}\label{phi'form}
 \phi'(t)=\sum_{j=1}^k \lambda_{i_j}'(t)=\sum_{j=1}^k \uu_{i_j}(t)^*B\uu_{i_j}(t).
 \end{equation}
 Since $\uu_1(t),\ldots,\uu_n(t)$ is an orthonormal basis in $\C^n$ Ky Fan inequality yields $\phi'(t)\le \sum_{j=1}^k \lambda_j(B)$.
 As $\lambda_1(t),\ldots,\lambda_n(t)$ are analytic in $\R\setminus Z$ and continuous on $\R$ it follows that for any real $t_0<t_1$
 \begin{equation}\label{integfor}
 \phi(t_1)-\phi(t_0)=\int_{t_0}^{t_1} \phi'(t')dt'\le (t_1-t_0)\sum_{j=1}^k \lambda_j(B).
 \end{equation}
 Choose $t_0=0, t_1=1$ to deduce the inequality \eqref{wielandt1}.  

 Assume that the equality \eqref{wieleq} holds.   Clearly, $(0,1)\setminus Z=\cup_{l=1}^N I_l \cap (0,1)$.  For simplicity of notation we let 
 $(0,1)\setminus Z=\cup_{l=1}^r (b_{l-1},b_l)$, where $b_0=0, b_r=1$.  Fix $\tau\in(b_{l-1}, b_l)$.
 Since $\phi(t)$ is analytic in $(b_{l-1},b_l)$ we deduce that
 \begin{equation}\label{condwieleq}
 \sum_{j=1}^k \uu_{i_j}(\tau)^* B\uu_{i_j}(\tau)=\sum_{j=1}^k \lambda_j(B).
 \end{equation} 
 Since $\phi(t)$ is continuous on $\R$ it follows that $\phi(t)=\sum_{j=1}^n (\lambda_{i_j}(A)+t\lambda_j(B))$ for $t\in[0,1]$.
 Let $\U(\tau)=\span(\uu_{i_1}(\tau),\ldots,\uu_{i_j}(\tau))$.  Ky Fan's theorem claims that $\U(\tau)$ is an invariant subspace of $B$ corresponding
 to the first $k$ eigenvalues of $B$.  Clearly $\U(\tau)$ is an invariant subspace of $A(\tau)$.  Hence $\U(\tau)$ is an invariant subspace $A$ of dimension
 $k$.   

 Let $\mu_{1,l}(t,\tau)\ge \cdots\ge \mu_{k,l}(t, \tau)$ be the eigenvalues of the restriction of $A(t)$ to $\U(\tau)$ for $t\in \R$.  
 Since $\lambda_1(t),\ldots,\lambda_n(t)$ are analytic in $(b_{l-1},b_l)$ it follows that $\mu_{i,l}(t,\tau)=\lambda_{p_j(\tau)}(t)$ for $t\in (b_{l-1},b_l)$
 and some integers $1\le p_1(\tau)<p_2(\tau)<\ldots<p_k(\tau)\le n$.  
 In view of \eqref{lambdtmultip} we can assume that $p_j(\tau)=i_j$ for $j\in [k]$ and each $\tau \in(b_{l-1},b_l)$.
 Fix $\tau\in(b_{l-1},b_l)$ and let $\U_l=\U(\tau)$.  Then $\mu_{j,l}(t):=\mu_{j}(t,\tau)$ is $\lambda_{i_j}(t)$ on $(b_{l-1},b_l)$ for $j\in[k]$.
 Since each $\mu_{j,l}(t)$ and $\lambda_{i_j}(t)$ are continuous on $[b_{l-1},b_l]$ we deduce that $\mu_{j,l}(t)=\lambda_{i_j}(t)$ for $j\in[k]$ and $t\in[b_{l-1},b_l]$.
 This shows that \eqref{wieleq} implies the existence of $\U_1,\ldots,\U_r$ with the claimed properties.

 Assume now $\U_1,\ldots, \U_r\subset \C^n$ are $k$-dimensional invariant subspaces of $A$ and $B$ satisfying the assumptions of the theorem.
 Let $A_l,B_l$ be the restricitons of $A$ and $B$ to $\U_l$.  Denote $A_l(t):=A_l+tB_l$.  Then 
 \[\tr A_l(t)=\sum_{j=1}^k \mu_{j,l}(t)=\tr A_l +t\tr B_l =\tr A(b_{l-1}) +(t-b_{l-1})\sum_{j=1}^k \lambda_j(B).\]
 Since $\tr A_l(t)=\sum_{j=1}^k \lambda_{i_j}(t)$ on each $[b_{l-1},b_l]$ we deduce \eqref{wieleq}.
 \qed

 \section{Remarks and an open problem}

 We remark that in the case of $k=1$ we can assume that $r=1$.  Indeed, on each $\U_l$ the pencil $A_l(t)$ has one analytic eigenvalue 
 $\lambda_{1,l} (t):=\lambda_{i(l)}+t\lambda_1(B)$.
 Hence for two distinct $l_1,l_2\in [r]$ either $\mu_{1,l_1}(t)$ is identically $\mu_{1,l_2}(t)$ or $\mu_{1,l_1}(t)-\mu_{1,l_2}(t)$ is never zero.
 We do not know if in a general case 
 we can always assume that $r=1$ in Theorem \ref{Wielineqeq}.  (For that one needs to discuss only the case where $r=2$.)

 Assume that \eqref{wieleq} holds.  Then there exists a $k$ dimensional subspace $\U_1$, invariant under $A$ and $B$
 such that the restriction of $\U_1$ is spanned by the first $k$ eigenvectors of $B$ and by the eigenvalues $1\le i_1<\cdots<i_k\le n$ of $A$. 
 The following proposition gives necessary and sufficient conditions on the existence of such $\U$:
 \begin{prop}\label{existU1}  Let $A,B\in\rH_n$ and let $1\le i_1<\cdots<i_k\le n$ be integers.  Then the following conditions are equivalent:
 \begin{enumerate}
 \item There exists a $k$-dimensional subspace $\U\subset\C^n$ satisfying the following properties:
 \begin{enumerate}
 \item $\U$ is an invariant subspace of $A$ spanned by $k$ eigenvectors corresponding 
 to the eigenvalues $\lambda_{i_1}(A),\ldots,\lambda_{i_k}(A)$.
 \item $\U$ is an invariant subspace of $B$ spanned by the eigenvectors corresponding to the first $k$ eigenvalues of $B$.
 \end{enumerate}
 \item There exist $t_1>0$ and $k$ integers $1\le p_1<\cdots<p_k\le n$ such that
 \begin{equation}\label{wieleqt1}
 \sum_{j=1}^k \lambda_{p_j}(A+t_1B)=\sum_{j=1}^k \lambda_{p_j}(A)+t_1\sum_{j=1}^k \lambda_j(B),\; \lambda_{p_j}(A)=\lambda_{i_j}(A), \;j\in[k].
 \end{equation}
  \item Let $\nu_1(A,B),\ldots,\nu_{n}(A,B)$ be defined as in Lemma \ref{pertform}.  Then there exist $k$ integers $1\le p_1<\cdots<p_k\le n$ such that
 \begin{equation}\label{nuABkid}
 \sum_{j=1}^k \nu_{p_j}(A,B)=\sum_{j=1}^k \lambda_j(B), \quad  \lambda_{p_j}(A)=\lambda_{i_j}(A), \;j\in[k].
 \end{equation}
 \end{enumerate}
 \end{prop}
 \proof  \emph{1.}$\Rightarrow$\emph{2.}  Let $A_1,B_1$ be a restriction of $A,B$ to $\U$ respectively.  Let $\mu_1(t)\ge \cdots\ge\mu_k(t)$ be the eigenvalues
 of $A_1(t)=A_1+t B_1$ for $t\in\R$.  Assume that $(0,t_1]\subset \R\setminus Z$ for some $t_1>0$.  Then $\mu_j(t)=\lambda_{p_j}(t)$ for $t\in (0,t_1]$ and $j\in[k]$.
 As $\U$ is an invariant subspace corresponding to $\lambda_{i_j}(A), j\in[k]$ and $\lambda_{p_j}(A), j\in [k]$
 we deduce the second part of  \eqref{wieleqt1}.  Since 
 \[\tr A_1(t)=\tr A_1+t\tr B_1=\sum_{j=1}^k \lambda_{i_j}(A)+t\sum_{j=1}^k \lambda_j(B)=\sum_{j=1}^k \lambda_{p_j}(A)+t\sum_{j=1}^k \lambda_j(B)\]
 we deduce the first part of \eqref{wieleqt1}.  

 \emph{2.}$\Rightarrow$\emph{3.}  Apply Theorem \ref{Wielineqeq} to $A$ and $t_1 B$.  Note that $\lambda_j(t_1B)=t_1\lambda_j(B)$ for $j\in[n]$.
 Let $\U=\U_1$.  Hence there exists $t_2\in(0,t_1)$ such that $(0,t_2)\in\R\setminus Z$.  Let $A_1,B_1,A_1(t)$ be as above.  Then $\mu_j(t)=\lambda_{p_j}(t)$
 for $j\in[k]$.  In particular, $\sum_{j=1}^k \lambda'_{p_j}(0^+)= \sum_{j=1}^k \nu_{p_j}(A,B)=\sum_{j=1}^k \lambda_j(B)$.  Use the second part of \eqref{wieleqt1}
 to deduce \eqref{nuABkid}.

 \emph{3.}$\Rightarrow$\emph{1.}   \eqref{derlambi0+} yields that 
 $\sum_{j=1}^k \lambda_{p_j}(0^+)=\sum_{j=1}^k \nu_{p_j}(A,B)$.  \eqref{difforeiglami+} implies that there exists an invariant subspace $\U$
 of $A$ spanned by orthonormal eigenvectors $\uu_{p_j}$ corresponding to the eigenvalue $\lambda_{p_j}(A)$ for $j\in[k]$ such that 
 $\sum_{j=1}^k \lambda_{p_j}(0^+)=\sum_{j=1}^k \uu_{p_j}^* B\uu_{p_j}$.  The first equality in \eqref{nuABkid} yields that $\U$ is an invariant subspace
 of $B$ corresponding to first $k$ eigenvalues of $B$.   The second equality in \eqref{nuABkid} yields \emph{1}.\qed 

 Observe that the condition \emph{3.} of Lemma \ref{existU1} can be verified efficiently.  Furthermore, given $\U$ which satisfies the condition \emph{1.} of Lemma 
 \ref{existU1} then there exist $1\le p_1<\cdots<p_k\le n$ and a maximal $t_1$, possibly $t_1=\infty$, such that the condition \emph{2.} holds for
 $t\in  [0,t_1]$ but not for $t>t_1$.

 Consider the following example.
 Assume that $\alpha_1>\alpha_2=\alpha_3, \beta_1>\beta_2\ge \beta_3$.  Suppose that $A=\diag(\alpha_1,\alpha_2,\alpha_3), B=\diag(\beta_3,\beta_1,\beta_2)$.
 Consider the case $k=1$ and $i_1=3$.  So $\U=\span((0,1,0)\trans)$.  Then  $p_1=2$ and $t_1=\frac{\alpha_1-\alpha_2}{\beta_1-\beta_3}$.

\end{document}